\documentclass[reqno]{amsart}
\usepackage{setspace,amssymb}
\usepackage{ifpdf}
\ifpdf
 \usepackage[hyperindex,pagebackref]{hyperref}%
\else
 \expandafter\ifx\csname dvipdfm\endcsname\relax
 \usepackage[hypertex,hyperindex,pagebackref]{hyperref}
 \else
 \usepackage[dvipdfm,hyperindex,pagebackref]{hyperref}
 \fi
\fi
\allowdisplaybreaks[4]
\theoremstyle{plain}
\newtheorem{thm}{Theorem}

\theoremstyle{remark}
\newtheorem{rem}{Remark}
\DeclareMathOperator{\td}{d}

\begin{document}

\title[An explicit formula for Bernoulli polynomials]
{An explicit formula for Bernoulli polynomials in terms of $\boldsymbol{r}$-Stirling numbers of the second kind}

\author[B.-N. Guo]{Bai-Ni Guo}
\address[Guo]{School of Mathematics and Informatics, Henan Polytechnic University, Jiaozuo City, Henan Province, 454010, China}
\email{\href{mailto: B.-N. Guo <bai.ni.guo@gmail.com>}{bai.ni.guo@gmail.com}, \href{mailto: B.-N. Guo <bai.ni.guo@hotmail.com>}{bai.ni.guo@hotmail.com}}
\urladdr{\url{https://www.researchgate.net/profile/Bai-Ni_Guo/}}

\author[I. Mez\H{o}]{Istv\'an Mez\H{o}}
\address[Mez\H{o}]{Department of Mathematics, Nanjing University of Information Science and Technology, Nanjing City, 210044, China}
\email{\href{mailto: I. Mez\H{o} <mezo.istvan@inf.unideb.hu>}{mezo.istvan@inf.unideb.hu}}
\urladdr{\url{http://www.inf.unideb.hu/valseg/dolgozok/mezoistvan/mezoistvan.html}}

\author[F. Qi]{Feng Qi}
\address[Qi]{Department of Mathematics, College of Science, Tianjin Polytechnic University, Tianjin City, 300387, China}
\email{\href{mailto: F. Qi <qifeng618@gmail.com>}{qifeng618@gmail.com}, \href{mailto: F. Qi <qifeng618@hotmail.com>}{qifeng618@hotmail.com}, \href{mailto: F. Qi <qifeng618@qq.com>}{qifeng618@qq.com}}
\urladdr{\url{http://qifeng618.wordpress.com}}

\begin{abstract}
In the paper, the authors establish an explicit formula for computing Bernoulli polynomials at non-negative integer points in terms of $r$-Stirling numbers of the second kind.
\end{abstract}



\thanks{This paper was typeset using \AmS-\LaTeX}

\maketitle

\section{Introduction}

It is well known that Bernoulli numbers $B_{k}$ for $k\ge0$ may be generated by
\begin{equation}\label{Bernumber-dfn}
\frac{x}{e^x-1}=\sum_{k=0}^\infty B_k\frac{x^k}{k!}=1-\frac{x}2+\sum_{k=1}^\infty B_{2k}\frac{x^{2k}}{(2k)!}, \quad |x|<2\pi
\end{equation}
and that Bernoulli polynomials $B_n(x)$ for $n\ge0$ and $x\in\mathbb{R}$ may be generated by
\begin{equation}\label{bp}
\frac{te^{xt}}{e^t-1}=\sum_{n=0}^\infty B_n(x)\frac{t^n}{n!}, \quad |t|<2\pi.
\end{equation}
\par
In combinatorics, Stirling numbers of the second kind $S(n,k)$ are equal to the number of partitions of the set $\{1,2,\dotsc,n\}$ into $k$ non-empty disjoint sets. Stirling numbers of the second kind $S(n,k)$ for $n\ge k\ge0$ may be computed by
\begin{equation}\label{Stirling-Number-dfn}
S(n,k)=\frac1{k!}\sum_{\ell=0}^k(-1)^{k-\ell}\binom{k}{\ell}\ell^{n}.
\end{equation}
In the paper~\cite{Broder}, among other things, Stirling numbers $S(n,k)$ were combinatorially generalized as $r$-Stirling numbers of the second kind, denoted by $S_r(n,k)$ here, for $r\in\mathbb{N}$, which may be alternatively defined as the number of partitions of the set $\{1,2,\dotsc,n\}$ into $k$ non-empty disjoint subsets such that the numbers $1,2,\dotsc,r$ are in distinct subsets.
\par
Note that
\begin{equation*}
S(0,0)=1,\quad S_0(n,k)=S(n,k),
\end{equation*}
and, when $n\in\mathbb{N}$,
\begin{equation*}
S(n,0)=0,\quad S_1(n,k)=S(n,k).
\end{equation*}
\par
In~\cite[p.~536]{GKP-Concrete-Math-1989} and~\cite[p.~560]{GKP-Concrete-Math-2nd}, the simple formula
\begin{equation}\label{Bernoulli-Stirling-eq}
B_n=\sum_{k=0}^n(-1)^k\frac{k!}{k+1}S(n,k), \quad n\in\{0\}\cup\mathbb{N}
\end{equation}
for computing Bernoulli numbers $B_n$ in terms of Stirling numbers of the second kind $S(n,k)$ was incidentally obtained.
Recently, four alternative proofs for the formula~\eqref{Bernoulli-Stirling-eq} were supplied in~\cite{ANLY-D-12-1238.tex} and its preprint~\cite{Bernoulli-Stirling2-3P.tex}. For more information on calculation of Bernoulli numbers $B_n$, please refer to the papers~\cite{exp-derivative-sum-Combined.tex, recursion, Bernoulli-No-Int-New.tex, Eight-Identy-More.tex, Tan-Cot-Bernulli-No.tex, CAM-D-13-01430-Xu-Cen}, especially to the article~\cite{Gould-MAA-1972}, and plenty of references therein.
\par
The aim of this paper is to generalize the formula~\eqref{Bernoulli-Stirling-eq}. Our main result may be formulated as the following theorem.

\begin{thm}\label{Bernoulli-Poly-r-Stirling-thm}
For all integers $n,r\ge0$, Bernoulli polynomials $B_n(r)$ may be computed in terms of $r$-Stirling numbers of the second kind $S_r(n+r,k+r)$ by
\begin{equation}\label{Bernoulli-Poly-r-Stirling-eq}
B_n(r)=\sum_{k=0}^n(-1)^k\frac{k!}{k+1}S_r(n+r,k+r).
\end{equation}
\end{thm}

In the final section of this paper, several remarks are listed.

\section{Proof of Theorem~\ref{Bernoulli-Poly-r-Stirling-thm}}

We are now in a position to verify our main result.
\par
For $n,r\ge0$, let
\begin{equation}\label{F(n-r)(x)-def}
F_{n,r}(x)=\sum_{k=0}^nk!S_r(n+r,k+r)x^k.
\end{equation}
By~\cite[p.~250, Theorem~16]{Broder}, we have
\begin{equation*}
\sum_{n=0}^\infty S_r(n+r,k+r)\frac{t^n}{n!}
=\sum_{n=k}^\infty S_r(n+r,k+r)\frac{t^n}{n!}
=\frac1{k!}e^{rt}(e^t-1)^k,
\end{equation*}
where $S_r(n,m)=0$ for $m>n$, see~\cite[p.~243, (10)]{Broder}. Accordingly, we obtain
\begin{align*}
\sum_{n=0}^\infty F_{n,r}(x)\frac{t^n}{n!}
&=\sum_{n=0}^\infty \sum_{k=0}^nk!x^kS_r(n+r,k+r)\frac{t^n}{n!}\\
&=\sum_{k=0}^\infty k!x^k\sum_{n=k}^\infty S_r(n+r,k+r)\frac{t^n}{n!}\\
&=e^{rt}\sum_{k=0}^\infty x^k(e^t-1)^k\\
&=\frac{e^{rt}}{1-x(e^t-1)}.
\end{align*}
Integrating with respect to $x\in[0,s]$ for $s\in\mathbb{R}$ on both sides of the above equation yields
\begin{equation}\label{tmp}
\sum_{n=0}^\infty \biggl[\int_0^s F_{n,r}(x)\td x\biggr] \frac{t^n}{n!}
=-e^{rt}\frac{\ln(1+s-se^t)}{e^t-1}.
\end{equation}
On the other hand,
\begin{equation*}
\int_0^s F_{n,r}(x)\td x=\sum_{k=0}^n\frac{k!}{k+1}S_r(n+r,k+r)s^{k+1}.
\end{equation*}
Substituting this into the equation~\eqref{tmp} concludes that
\begin{equation*}
\sum_{n=0}^\infty \sum_{k=0}^n\frac{k!}{k+1}S_r(n+r,k+r)s^{k+1} \frac{t^n}{n!}
=-e^{rt}\frac{\ln(1+s-se^t)}{e^t-1}.
\end{equation*}
Taking $s=-1$ in the above equation and making use of the generating function~\eqref{bp} result in
\begin{equation*}
\sum_{n=0}^\infty \Biggl[\sum_{k=0}^n(-1)^{k+1} \frac{k!}{k+1}S_r(n+r,k+r)\Biggr]\frac{t^n}{n!}
=-\frac{te^{rt}}{e^t-1}
=\sum_{n=0}^\infty [-B_n(r)]\frac{t^n}{n!},
\end{equation*}
which implies the formula~\eqref{Bernoulli-Poly-r-Stirling-eq}. The proof of Theorem~\ref{Bernoulli-Poly-r-Stirling-thm} is complete.

\section{Remarks}

Finally we would like to give several remarks on Theorem~\ref{Bernoulli-Poly-r-Stirling-thm} and its proof.

\begin{rem}
Since $B_n(0)=B_n$ and $S_0(n,k)=S(n,k)$, when $r=0$, the formula~\eqref{Bernoulli-Poly-r-Stirling-eq} becomes~\eqref{Bernoulli-Stirling-eq}. Therefore, our Theorem~\ref{Bernoulli-Poly-r-Stirling-thm} generalizes the formula~\eqref{Bernoulli-Stirling-eq}.
\end{rem}

\begin{rem}
It is easy to see that
\begin{equation*}
F_{n,0}(1)=\sum_{k=0}^nk!S(n,k),
\end{equation*}
which are just the classical ordered Bell numbers. For more information, please refer to the papers~\cite{Can-Joyce-JCTA-2012, Bell-Stirling-Lah-simp.tex} and closely related references therein.
\end{rem}

\begin{rem}
In the PhD thesis~\cite{MezoThesis}, the second author defined a variant of the polynomials $F_{n,r}(x)$. Hence, a simple combinatorial study and interpretation of the polynomials $F_{n,r}(x)$ is available therein.
\end{rem}

\end{document}